\documentclass{article}
\pdfoutput=1

\usepackage{arxiv}

\usepackage[utf8]{inputenc} % allow utf-8 input
\usepackage[english]{babel}
\usepackage[T1]{fontenc}    % use 8-bit T1 fonts
\usepackage{hyperref}       % hyperlinks
\usepackage{url}            % simple URL typesetting
\usepackage{booktabs}       % professional-quality tables
\usepackage{amsfonts}       % blackboard math symbols
\usepackage{nicefrac}       % compact symbols for 1/2, etc.
\usepackage{microtype}      % microtypography
\usepackage{amsmath}
\usepackage{amssymb}

\title{Explicit Multi-point Taylor Polynomial}
\author{
 Andrés Gómez Arias \\
  Facultad de Ciencias\\
  Universidad Nacional Autónoma de México\\
  \texttt{andresgz@ciencias.unam.mx} \\
}
\date{\today}

\begin{document}

\maketitle

\begin{abstract}
    The multi-point Taylor polynomial, which is the general, unique and of minimum degree ($mk+m-1$) polynomial $P_{k,m}(x)$ which interpolates a function's derivatives in multiple points is presented in its explicit form. A proof that this expression satisfies the multi-point Taylor polynomial's defining property is given. Namely, it is proven that for a k-differentiable function $f$ and a set of different m-points $\{a_1,...,a_m\}$, this polynomial satisfies $P^{(n)}_{k,m}(a_i) = f^{(n)}(a_i)\quad \forall \, i = 1,...,m\quad \&\quad \forall \, n = 0,...,k$. A discussion regarding previous expressions presented in the literature, which mostly consisted in recursion formulas and not explicit formulas, is made.
\end{abstract}


\section{Introduction}

It is well known that an infinitely-differentiable function $f:I\longrightarrow \mathbb{R}$ can represented by its Taylor series expansion around a single point $a\in I$. This in turn defines the single-point Taylor polynomial (STP) as the partial sum of the series, given by
\begin{equation}
\label{Taylor}
    p_k(x) = \sum_{n=0}^{k}\frac{(x-a)^{n}}{n!}f^{(n)}(a),
\end{equation}
which produces a good approximation of the function for increasing values of $k$ and in fact exactly interpolates the function's derivatives at the point $a$. That is, the STP is the unique and of minimal (of order $k$) polynomial with the property that 
\begin{equation}
    p^{(n)}_k(a) = f^{(n)}(a) \quad \forall\, n = 0,...,k.
\end{equation}
Furthermore, the STP can be defined for functions that are just $k$-differentiable, which together with the previous properties makes it very suitable for polynomial interpolation.

 On the other hand, for a set of different points $\mathcal{A} = \{a_1,...,a_m\}\subset I$, the Lagrange polynomial (LP)
\begin{equation}
\label{Lagrange}
    q_{m}(x) = \sum_{g=1}^{m}\prod^{m}_{\substack{h=1\\
    h\neq g}}\frac{(x-a_h)}{(a_g-a_h)}f(a_g),
\end{equation}
gives the unique polynomial (of minimum order $m$)  that interpolates the value of the function at the given points
\begin{equation}
    q_m(a_i) = f(a_i) \quad \forall\, i = 1,...,m,
\end{equation}
and can be found as a basic multi-point interpolating polynomial in many numerical analysis textbooks~\cite{Richard,Philip,Bulirsch,Hildebrand}.

One can go further and ask for a polynomial which combines both of these properties, defining the multi-point Taylor polynomial (MTP). That is, for a k-differentiable real function $f:I\longrightarrow \mathbb{R}$ and a set of different points $\mathcal{A} = \{a_1,...,a_m\}\subset I$, the MTP is the unique and minimal polynomial (of order $mk+m-1$) which fulfills the conditions
\begin{equation}
\label{property}
    P^{(n)}_{k,m}(a_i) = f^{(n)}(a_i)\quad \forall \, i = 1,...,m\quad \&\quad \forall \, n = 0,...,k.
\end{equation}
The solution to these conditions is also known as the Hermite or Osculatory interpolation polynomial, and multiple standard methods and iterative algorithms can be found for producing it~\cite{Richard,Philip,Bulirsch,Note_2Point}, but these do not give its explicit and final form. This is what will be presented on this paper.

The use of MTP or similar interpolants has encountered many applications. For example, the authors in Ref.~\cite{NASA} presented a two-point Taylor series expansion whose coefficients can be iterated to produce higher order expressions and applied this to the two-body problem, which consisted of an expansion of the solutions both at the perigee and apogee. Their results turned out to be significantly better than a basic Taylor approximation around just the perigee. Another study~\cite{Lopez} used a convenient form of the MTP and applied it to the approximation of solutions of second order linear differential equations on an interval with boundary conditions on the extremes, which was done by approximating the solutions around the two boundary points. On the other hand, the study~\cite{Diff_eq} performed a multiple point expansion by requesting the satisfaction of the boundary conditions and of the differential equation at a cloud of points inside the domain. Finally, the work in~\cite{Image_comp} used a multi-point Taylor series formula in terms of general basis functions in order to construct an image compression/decompression method.

\section{The multi-point Taylor polynomial}

The explicit form of the MTP will now presented in the form of a theorem and is the main result of this paper.

\textbf{Theorem:}

For a k-differentiable real function $f:I\longrightarrow \mathbb{R}$ and a set of different points $\mathcal{A} = \{a_1,...,a_m\}\subset I$, the explicit expression for the MTP is
\begin{equation}
\label{Polynomial}
    P_{k,m}(x) = \sum_{g=1}^{m}\bigg[\prod^{m}_{\substack{h=1\\
    h\neq g}}\frac{(x-a_h)}{(a_g-a_h)}\bigg]^{k+1}\sum_{n=0}^{k}\frac{(x-a_g)^{n}}{n!}F^{n,g}_{k,m}[\mathcal{A}],
\end{equation}
where $F^{n,g}_{k,m}[\mathcal{A}]$ are constant factors given by
\begin{equation}
    F^{n,g}_{k,m}[\mathcal{A}] = \sum_{j_1+...+j_m = n}\binom{n}{j_1,...,j_m}f^{(j_g)}(a_g)\prod^{m}_{\substack{l=1\\
    l\neq g}}\frac{(k+j_l)!}{k!}(a_l-a_g)^{-j_l}.
\end{equation}
This polynomial satisfies the conditions (\ref{property}) and thus interpolates a function's derivatives at multiple points.

Under the understanding that empty sums are null and empty products are one, this expression recovers the STP of Eq. (\ref{Taylor}) for $m = 1$, while $k=0$  recovers the LP of Eq. (\ref{Lagrange}). This is a straightforward property which comes from the fact that the conditions (\ref{property}) recover the STP's or the LP's conditions respectively for each case.

\textbf{Proof:}

It will now be proven that the polynomial $P_{k,m}(x)$ of Eq. (\ref{Polynomial}) satisfies the conditions (\ref{property}). First notice that $P^{(n)}_{k,m}(a_i)$ will only have contributions from the $g=i$ term of the first sum, as the order of the derivative is not big enough to suppress the $(x-a_i)$ factor involved in the other terms which involve
\begin{equation}
    \bigg[\prod^{m}_{\substack{h=1\\
    h\neq g}}\frac{(x-a_h)}{(a_g-a_h)}\bigg]^{k+1},
\end{equation}
with $g\neq i$, as it is assumed that $n<k+1$. Therefore,

\begin{equation}
\label{a}
\begin{split}
    P^{(n)}_{k,m}(a_i) = \frac{\mathrm{d}^{n}}{\mathrm{d}x^{n}}\Bigg[\bigg[\prod^{m}_{\substack{h=1\\
    h\neq i}}\frac{(x-a_h)}{(a_i-a_h)}\bigg]^{k+1}\sum_{s=0}^{k}&\frac{(x-a_i)^{s}}{s!}F^{s,i}_{k,m}[\mathcal{A}]\Bigg]_{x=a_i}. \\
\end{split}
\end{equation}

With that, using the well known Leibniz rule for higher order derivatives
\begin{equation}
    (fg)^{(n)} = \sum_{r=0}^{n}\binom{n}{r}f^{(n-r)}g^{(r)},
\end{equation}
then
\begin{equation}
\label{b}
\begin{split}
\frac{\mathrm{d}^{n}}{\mathrm{d}x^{n}}\Bigg[\bigg[\prod^{m}_{\substack{h=1\\
    h\neq i}}\frac{(x-a_h)}{(a_i-a_h)}\bigg]^{k+1}\sum_{s=0}^{k}\frac{(x-a_i)^{s}}{s!}A_s\Bigg]_{x=a_i} & =  \sum_{r=0}^{n}\binom{n}{r}\frac{\mathrm{d}^{r}}{\mathrm{d}x^{r}}\Bigg[\bigg[\prod^{m}_{\substack{h=1\\
    h\neq i}}\frac{(x-a_h)}{(a_i-a_h)}\bigg]^{k+1}\Bigg]_{x=a_i} \\
    &\hspace{4cm}\times \sum_{s=0}^{k}\frac{(a_i-a_i)^{s-n+r}}{(s-n+r)!}F^{s,i}_{k,m}[\mathcal{A}]. \\
\end{split}
\end{equation}

Again, only the $s = n-r$ term survives. Now, the general Leibniz rule for more than two factors is
\begin{equation}
    \bigg(\prod_{h=1}^{m} f_h\bigg)^{(n)} = \sum_{k_1+...+k_m = n}\binom{n}{k_1,...,k_m}\prod_{h=1}^{m} f_h^{(k_h)}
\end{equation}
(a proof for this relation can be found for example in Ref. \cite{Gen_Leibniz}), so
\begin{equation}
\begin{split}
    \frac{\mathrm{d}^{r}}{\mathrm{d}x^{r}}\Bigg[\prod^{m}_{\substack{h=1\\
    h\neq i}}\frac{(x-a_h)^{k+1}}{(a_i-a_h)^{k+1}}\Bigg]_{x=a_i} & = \sum_{k_1+...+k_m = r}^{i}\binom{r}{k_1,...,k_m}^i\prod_{\substack{h=1\\
    h\neq i}}^{m} \frac{(k+1)!}{(k+1-k_h)!}\frac{(a_i-a_h)^{k+1-k_h}}{(a_i-a_h)^{k+1}} \\
    & = [(k+1)!]^{m-1}\sum_{k_1+...+k_m = r}^{i}\binom{r}{k_1,...,k_m}^i\prod_{\substack{h=1\\
    h\neq i}}^{m} \frac{1}{(k+1-k_h)!}(a_i-a_h)^{-k_h},
\end{split}
\end{equation}
where the multinomial sum and multinomial coefficient with upper index $\sum_{k_1+...+k_m = r}^{i}$ and $\binom{r}{k_1,...,k_m}^i$ denote that $k_i$ is not taken into account. Joining this with equations (\ref{a}) and (\ref{b}), it follows that

\begin{equation}
\begin{split}
P^{(n)}_{k,m}(a_i) & =  \sum_{r=0}^{n}\binom{n}{r} [(k+1)!]^{m-1}F^{n-r,i}_{k,m}[\mathcal{A}]\sum_{k_1+...+k_m = r}^{i}\binom{r}{k_1,...,k_m}^i\prod_{\substack{h=1\\
    h\neq i}}^{m} \frac{1}{(k+1-k_h)!}(a_i-a_h)^{-k_h}. \\
\end{split}
\end{equation}
Now,
\begin{equation}
\begin{split}
    \binom{n}{r}[(k+1)!]^{m-1}\binom{r}{k_1,...,k_m}^i&F^{n-r,i}_{k,m}[\mathcal{A}] \\
    & = \frac{(k+1)^{m-1}}{(n-r)!}\binom{n}{k_1,...,k_m}^i\sum_{j_1+...+j_m = n-r}\binom{n-r}{j_1,...,j_m}f^{(j_i)}(a_i)\prod^{m}_{\substack{l=1\\
    l\neq i}}\frac{(k+j_l)!}{(a_l-a_i)^{j_l}} \\
    & = (k+1)^{m-1}\binom{n}{k_1,...,k_m}^i\sum_{j_1+...+j_m = n-r}\binom{1}{j_1,...,j_m}f^{(j_i)}(a_i)\prod^{m}_{\substack{l=1\\
    l\neq i}}\frac{(k+j_l)!}{(a_l-a_i)^{j_l}}. \\
\end{split}
\end{equation}
Therefore, the evaluation of the polynomial can be expressed as the sum
\begin{equation}
\label{c}
\begin{split}
P^{(n)}_{k,m}(a_i) & =  (k+1)^{m-1}\sum_{r=0}^{n}B_r, \\
\end{split}
\end{equation}
with
\begin{equation}
    B_r = \sum_{k_1+...+k_m = r}^{i}\binom{n}{k_1,...,k_m}^i\sum_{j_1+...+j_m = n-r}\binom{1}{j_1,...,j_m}f^{(j_i)}(a_i)\prod^{m}_{\substack{l=1\\
    l\neq i}}\prod_{\substack{h=1\\
    h\neq i}}^{m}\frac{(k+j_l)!}{(k+1-k_h)!(a_i-a_h)^{k_h}(a_l-a_i)^{j_l}}.
\end{equation}

As the products range between the same indices, these can be joined together under the same label. Furthermore, when grouping the last two denominators, a global sign appears, which is given by
\begin{equation}
    \prod_{\substack{h=1\\
    h\neq i}}^{m}(-1)^{k_h} = (-1)^{r},
\end{equation}
because $k_1+...+k_m = r$, where $k_i$ is not considered. The result of this process is
\begin{equation}
    B_r = (-1)^{r}\sum_{k_1+...+k_m = r}^{i}\binom{n}{k_1,...,k_m}^i\sum_{j_1+...+j_m = n-r}\binom{1}{j_1,...,j_m}f^{(j_i)}(a_i)\prod_{\substack{h=1\\
    h\neq i}}^{m}\frac{(k+j_h)!}{(k+1-k_h)!}(a_h-a_i)^{-j_h-k_h}.
\end{equation}

The only term that could produce the case $j_i = n$, so that the expression includes a $f^{(n)}(a_i)$ term, is the one where $r=0$. Making an explicit calculation of this term,
\begin{equation}
\begin{split}
    B_0 & = (-1)^{0}\sum_{k_1+...+k_m = 0}^{i}\binom{n}{k_1,...,k_m}^i\sum_{j_1+...+j_m = n}\binom{1}{j_1,...,j_m}f^{(j_i)}(a_i)\prod_{\substack{h=1\\
    h\neq i}}^{m}\frac{(k+j_h)!}{(k+1-k_h)!}(a_h-a_i)^{-j_h-k_h} \\
    & = \sum_{j_1+...+j_m = n}\binom{n}{j_1,...,j_m}f^{(j_i)}(a_i)\prod_{\substack{h=1\\
    h\neq i}}^{m}\frac{(k+j_h)!}{(k+1)!}(a_h-a_i)^{-j_h}. \\
\end{split}
\end{equation}
Separating the sought $j_i = n$ term from the rest,
\begin{equation}
\begin{split}
    B_0 & = f^{(n)}(a_i)\prod_{\substack{h=1\\
    h\neq i}}^{m}\frac{k!}{(k+1)!}(a_h-a_i)^{0} + \sum_{j_1+...+j_m = n}^{j_i\neq n}\binom{n}{j_1,...,j_m}f^{(j_i)}(a_i)\prod_{\substack{h=1\\
    h\neq i}}^{m}\frac{(k+j_h)!}{(k+1)!}(a_h-a_i)^{-j_h} \\
    & = \frac{1}{(k+1)^{m-1}}f^{(n)}(a_i)+ \sum_{j_1+...+j_m = n}^{j_i\neq n}\binom{n}{j_1,...,j_m}f^{(j_i)}(a_i)\prod_{\substack{h=1\\
    h\neq i}}^{m}\frac{(k+j_h)!}{(k+1)!}(a_h-a_i)^{-j_h}, \\
\end{split}
\end{equation}
which gives the right term when inserted in equation (\ref{c}). Thus, now it must only be proven that all the other terms cancel each other. For this, it is convenient to group every coefficient of the derivatives of $f(a_i)$. That is,
\begin{equation}
\label{Expansion}
\begin{split}
    P^{(n)}_{k,m}(a_i) & =  (k+1)^{m-1}\sum_{r=0}^{n}B_r =  f^{(n)}(a_i)+n!\sum_{j=0}^{n-1}\frac{f^{(j)}(a_i)}{(j)!}C_{j},
\end{split}
\end{equation}
where, for $0\leq j \leq n-1$,
\begin{equation}
\begin{split}
    C_{j} & = \frac{(k+1)^{m-1}}{n!}\sum_{r=0}^{n-j}(-1)^{r}\sum_{k_1+...+k_m = r}^{i}\binom{n}{k_1,...,k_m} \\
    &\hspace{6cm}\times\sum_{j_1+...+j_m = n-r-j}^{i}\binom{1}{j_1,...,j_m}\prod_{\substack{h=1\\
    h\neq i}}^{m}\frac{(k+j_h)!}{(k+1-k_h)!}(a_h-a_i)^{-j_h-k_h} \\
    & =(k+1)^{m-1}\sum_{r=0}^{n-j}(-1)^{r}\sum_{\substack{k_1+...+k_m = r \\
    j_1+...+j_m = n-r-j}}^{i}\prod_{\substack{h=1\\
    h\neq i}}^{m}\frac{(k+j_h)!}{(k+1-k_h)!k_h!j_h!}(a_h-a_i)^{-j_h-k_h} \\
    & = \sum_{r=0}^{n-j}(-1)^{r}\sum_{\substack{k_1+...+k_m = r \\
    j_1+...+j_m = n-r-j}}^{i}\prod_{\substack{h=1\\
    h\neq i}}^{m}\binom{k+j_h}{j_h}\binom{k+1}{k_h}(a_h-a_i)^{-j_h-k_h}. \\
\end{split}
\end{equation}

These terms automatically include the $j_i\neq n$ term from $B_0$, whose nullity was left pendent to be proven. The restricted limit on the sum over $r$ comes from the fact that $j> n-r$ is not taken. The notation $\sum_{j_1+...+j_m = n-r-j}^{i}$ is used again to denote that $j_i$ is not considered (in fact $j_i = j$). This applies for the sum over $k_i$ as well. The task is now to prove that each factor $C_j$ is equal to zero.

Now, because
\begin{equation}
    \sum_{r=0}^{n-j}(-1)^{r}\sum_{\substack{k_1+...+k_m = r \\
    j_1+...+j_m = n-r-j}}^{i}\big[...\big] = \sum_{S(k_h)+S(j_h) = n-j}^{i}(-1)^{S(k_h)}\big[...\big],
\end{equation}
where $S(k_h) = k_1+...+k_m$ and $S(j_h) = j_1+...+j_m$ (without $k_i$ nor $j_i$), it follows that
\begin{equation}
\begin{split}
    C_{j} & = \sum_{S(k_h)+S(j_h) =n_j}^{i}(-1)^{S(k_h)}\prod_{\substack{h=1\\
    h\neq i}}^{m}\binom{k+j_h}{j_h}\binom{k+1}{k_h}(a_h-a_i)^{-j_h-k_h}, \\
\end{split}
\end{equation}
where the positive integer $n_j = n-j$ was defined. It is now convenient to further join together the terms of the same power over $(a_h-a_i)$. That is, the terms with $j_w+k_w = s_w$, for $w = 1,...,m$ and $w\neq i$ will now be grouped, which is done by performing the sum over the possible $s_w$ terms. This results in
\begin{equation}
\begin{split}
    C_{j} & = \sum_{s_1+...+s_m = n_j}^{i} \sum_{\substack{j_v+k_v=s_v\\
    v=1,...,m}}^{i}\prod_{\substack{h=1\\
    h\neq i}}^{m}(-1)^{k_h}\binom{k+j_h}{j_h}\binom{k+1}{k_h}(a_h-a_i)^{-j_h-k_h} \\
    & = \sum_{s_1+...+s_m = n_j}^{i}\bigg[\prod_{\substack{w=1\\
    w\neq i}}^{m}(a_w-a_i)^{-s_w}\bigg] \sum_{\substack{j_v+k_v=s_v\\
    v=1,...,m}}^{i}\prod_{\substack{h=1\\
    h\neq i}}^{m}(-1)^{k_h}\binom{k+j_h}{j_h}\binom{k+1}{k_h}. \\
\end{split}
\end{equation}
Once again, the upper index on the sum denotes that the respective term with sub-index equal to $i$ is not considered. Having grouped every non-compatible term, it is now left to be proven that the factor of each one is zero. That is, one must prove that, for any fixed combination $\mathcal{S} = \{s_1,...,s_m\}$ such that $s_1+...+s_m = n_j$ (without considering $s_i$), the term
\begin{equation}
    D_{n_j,k}(\mathcal{S}) \equiv \sum_{\substack{j_v+k_v=s_v\\
    v=1,...,m}}^{i}\prod_{\substack{h=1\\
    h\neq i}}^{m}(-1)^{k_h}\binom{k+j_h}{j_h}\binom{k+1}{k_h}
\end{equation}
is identically zero, where $1 \leq n_j \leq n \leq k$. This will be done by strong induction over $n_j$. For $n_j=1$ the only possible combination is when only one element of $\mathcal{S}$ is equal to one, while the others are zero. Without loss of generality, this element can be taken to be $s_1$. Then, the sum
\begin{equation}
    \sum_{\substack{j_v+k_v=s_v\\
    v=1,...,m}}^{i}
\end{equation}
has two terms: One where $(j_1,k_1)=(1,0)$ and another where $(j_1,k_1)=(0,1)$, where in both cases all of the other $(j_v,k_v)$ pairs are zero. Thus
\begin{equation}
\begin{split}
    D_{1,k}(\mathcal{S}) & = \bigg[(-1)^{0}\binom{k+1}{1}\binom{k+1}{0}+(-1)^{1}\binom{k+0}{0}\binom{k+1}{1}\bigg]\prod_{\substack{h=2\\
    h\neq i}}^{m}(-1)^{0}\binom{k+0}{0}\binom{k+1}{0} \\
    & = (k+1)-(k+1) \\
    & = 0.
\end{split}
\end{equation}
Now the inductive step will be performed. Suppose that, for any integer $n'$ such that $1\leq n'\leq n_j$, any combination $\mathcal{S} = \{s_1,...,s_m\}$ such that $s_1+...+s_m = n'$ (without considering $s_i$) produces a $D_{n',k}(\mathcal{S})$ term which vanishes. Take $n_j+1$ and an arbitrary combination $\mathcal{S}' = \{s'_1,...,s'_m\}$ such that $s'_1+...+s'_m = n_j+1$ (without considering $s'_i$). Now it must be proven that $D_{n_j+1,k}(\mathcal{S}')$ is zero.

It is clear that there exists a combination $\mathcal{S} = \{s_1,...,s_m\}$ such that $s_1+...+s_m = n_j$ and all of the elements of $\mathcal{S}$ are equal to the ones of $\mathcal{S}'$ with the exception of one, which differs by one unity. This can be done by just taking a copy of $\mathcal{S}'$ except for one element $s'$ greater than zero (which always exists as $n'\geq 1$), for which the element $s'-1$ is taken instead. Without loss of generality, this element can be taken as the first one. Thus,
\begin{equation}
    \mathcal{S}' = \{s_1+1,s_2,...,s_m\},
\end{equation}
and therefore
\begin{equation}
\begin{split}
    D_{n_j+1,k}(\mathcal{S}') & = \sum_{j_1+k_1=s_1+1}\sum_{\substack{j_v+k_v=s_v\\
    v=2,...,m}}^{i}\prod_{\substack{h=1\\
    h\neq i}}^{m}(-1)^{k_h}\binom{k+j_h}{j_h}\binom{k+1}{k_h} \\
    & = \sum_{j_1+k_1=s_1+1}(-1)^{k_1}\binom{k+j_1}{j_1}\binom{k+1}{k_1}\sum_{\substack{j_v+k_v=s_v\\
    v=2,...,m}}^{i}\prod_{\substack{h=2\\
    h\neq i}}^{m}(-1)^{k_h}\binom{k+j_h}{j_h}\binom{k+1}{k_h}. \\
\end{split}
\end{equation}
The second sum of this equation can be expressed as a particular case of $D$. That is, if one takes the combination $\mathcal{S}'' = \{0,s_2...,s_m\}$, then $0+s_2+...+s_m = n_j-s_1$ and thus
\begin{equation}
\begin{split}
    D_{n_j-s_1,k}(\mathcal{S}'') & = (-1)^0\binom{k+0}{0}\binom{k+1}{0}\sum_{\substack{j_v+k_v=s_v\\
    v=2,...,m}}^{i}\prod_{\substack{h=2\\
    h\neq i}}^{m}(-1)^{k_h}\binom{k+j_h}{j_h}\binom{k+1}{k_h} \\
    & = \sum_{\substack{j_v+k_v=s_v\\
    v=2,...,m}}^{i}\prod_{\substack{h=2\\
    h\neq i}}^{m}(-1)^{k_h}\binom{k+j_h}{j_h}\binom{k+1}{k_h}. \\
\end{split}
\end{equation}
This is because the sum over $(j_1,k_1)$ has only one term, where the values of these variables are zero, and thus can be factored out. This implies that
\begin{equation}
\begin{split}
    D_{n_j+1,k}(\mathcal{S}') & = \sum_{j_1+k_1=s_1+1}(-1)^{k_1}\binom{k+j_1}{j_1}\binom{k+1}{k_1}D_{n_j-s_1,k}(\mathcal{S}''). \\
\end{split}
\end{equation}
There are two distinct possibilities. If $s_1<n_j$, then $n_j-s_1 \geq 1$ and thus $D_{n_j-s_1,k}(\mathcal{S}'')$ fulfills the requirements for the inductive step. This means that $D_{n_j-s_1,k}(\mathcal{S}'')$ vanishes and so $D_{n_j+1,k}(\mathcal{S}')$ vanishes too. The other case, for which $s_1 = n_j$, does not fulfills the inductive step conditions (one needs $n_j>0$) and in fact turns out to be equal to one: $D_{0,k}(\mathcal{S}'') = 1$. This in turn implies that
\begin{equation}
\label{Non_trivial}
\begin{split}
    D_{n_j+1,k}(\mathcal{S}') & = \sum_{j_1+k_1=n_j+1}(-1)^{k_1}\binom{k+j_1}{j_1}\binom{k+1}{k_1} \\
    & = \sum_{k_1=0}^{n_j+1}(-1)^{k_1}\binom{k+n_j+1-k_1}{n_j+1-k_1}\binom{k+1}{k_1}. \\
\end{split}
\end{equation}
Although non-trivial, this expression is more manageable and its nullity can be proven once again by induction. This is done in \hyperref[AppendixA]{Appendix A}. With this, in either case, $D_{n_j+1,k}(\mathcal{S}')$ vanishes and the inductive step is fulfilled.

Having proven that $D_{n_j,k}(\mathcal{S})$ vanishes for any integer $n_j$ and any combination $\mathcal{S}$ with the given conditions, this consequently proves that $C_j = 0$ for $j = 0,...,n-1$ and, by means of Eq. (\ref{Expansion}), finalizes the proof that (\ref{Polynomial}) fulfills the conditions (\ref{property}) for the MTP. $\blacksquare$

\section{Discussion and conclusions}
A benefit of using a MTP instead of a STP is that one can obtain a greater accuracy for the interpolation of a function with the knowledge of less derivatives. That is, if one knows the first $k$ derivatives of a function at $m$ different points, the MTP possesses the same degree (and thus the same accuracy) with respect to the STP where a higher number of derivatives ($mk+m-1$) is known.

Although generally expensive to be computed, there can be parts of the expression (\ref{Polynomial}) which can be calculated without knowledge of the function to be interpolated, and so can be obtained beforehand and saved for optimization. Such terms include
\begin{equation}
\label{comp1}
    \bigg[\prod^{m}_{\substack{h=1\\
    h\neq g}}\frac{(x-a_h)}{(a_g-a_h)}\bigg]^{k+1}
\end{equation}
(which are powers of LP's factors), and
\begin{equation}
\label{comp2}
    \prod^{m}_{\substack{l=1\\
    l\neq i}}\frac{(k+j_l)!}{k!}(a_l-a_i)^{-j_l},
\end{equation}
where the second one is to be calculated for $0\leq j_l \leq k$.

As for the limiting behaviour, this is explored in Ref.~\cite{franssens1999} for general interpolating basis functions, where convergence to the function is found in appropriate analytic regions of the function (i) when the number of derivatives $k$ tends to infinity for fixed $m$, and (ii) when the spacing between the points becomes zero for fixed $k$.

Apart from discussing convergence of the series with $k$ going to infinity, the work on Ref.~\cite{fine1962} gives a general method for obtaining the MTP for analytic functions in terms of an expansion
\begin{equation}
    \sum_{t=0}^{k}[p(z)]^{t}p_{m-1,t}(z),
\end{equation}
where
\begin{equation}
    p(z) = \prod_{h=1}^{m}(z-a_h),
\end{equation}
and $p_{m-1,t}(z)$ are polynomials of degree $m-1$ whose coefficients can be calculated in terms of a contour integral that could be generally computed via the Cauchy residue theorem. Nonetheless, an explicit expression for such coefficients is not given.

An explicit expression for the MTP was in fact already introduced by J. L. López and N. M. Temme for analytical functions, first for two distinct points~\cite{2Point_N_Temme} and later for an arbitrary number of points (some of which could also be repeated)~\cite{MPoint_N_Temme}. In this work, the authors presented a Taylor series expansion at different points, considering the remainder of the finite polynomial and convergence radius. In the context of this paper their result applied for a finite polynomial and for different points is
\begin{equation}
\label{Q_km}
    P_{k,m}(z) = \sum_{n=0}^{k}q_{n,m}(z)\prod_{l=1}^{m}(z-a_l)^{n},
\end{equation}
where $q_{n,m}(z)$ are polynomials of degree $m-1$ given by
\begin{equation}
    q_{n,m}(z) \equiv \sum_{g=1}^{m}A_{n,g}\prod^{m}_{\substack{h=1\\
    h\neq g}}\frac{(z-a_h)}{(a_g-a_h)},
\end{equation}
and the constants $A_{n,g}$ can be given by the Cauchy integral
\begin{equation}
\label{A_ng_Cauchy}
    A_{n,g} = \frac{1}{2\pi i}\int_{\mathcal{C}}\frac{f(w)\,\mathrm{d}w}{(w-a_g)\prod_{h=1}^{m}(w-a_h)^n},
\end{equation}
where the contour of integration $\mathcal{C}$ is a simple closed loop which encircles all of the points of the set $\mathcal{A}$ of points in the counterclockwise direction and is contained in the analytical region $\Omega$ of the function. The aforementioned explicit expression for the MTP is obtained when one takes an alternate expression for these constants, given by the composite derivative
\begin{equation}
\label{A_ng2}
    A_{n,g} = \bigg[\frac{1}{n!}\frac{\mathrm{d}^{n}}{\mathrm{d}w^{n}}\frac{f(w)}{\prod^{m}_{\substack{h=1\\
    h\neq g}}(w-a_h)^{n}}\bigg]\bigg|_{w=a_g}+\sum_{\substack{l=1\\
    l\neq g}}^{m}\bigg[\frac{1}{(n-1)!}\frac{\mathrm{d}^{n-1}}{\mathrm{d}w^{n-1}}\frac{f(w)}{(w-a_g)\prod^{m}_{\substack{h=1\\
    h\neq l}}(w-a_h)^{n}}\bigg]\bigg|_{w=a_l},
\end{equation}
and which was also presented in this work. This formula holds for $n=1,2,3...$, while for the $n=0$ case one can directly obtain that $A_{0,g} = a_g$ as Eq. (\ref{A_ng_Cauchy}) results in the trivial Cauchy formula. Furthermore, the study also presents the multi-point generalization of the disk as the region of convergence for one point, given by
\begin{equation}
    O_m \equiv \{ z\in \Omega, \, \prod_{h=1}^{m}|z-a_h| < r\}, \qquad \qquad r \equiv \mathrm{Inf}_{w \in \mathbb{C}/\Omega}\bigg\{\prod_{h=1}^{m}|w-a_h|\bigg\}.
\end{equation}
It is also worth noting that a general Laurent series for multiple points is also analyzed, which is something that was not considered here.

Although the explicit expression of the polynomial (\ref{Q_km}) is not similar to the one presented here, the uniqueness of the solution must ensure that, under certain algebraic manipulations, Eqs. (\ref{Polynomial}) and (\ref{Q_km}) must be equal. The main difference between both results is that the expansion of both polynomials is performed in a different manner, which produces different expressions for the coefficients of each term. That is, the constants $A_{n,g}$ are expressed as n-th derivatives of quotients between the function $f(w)$ and powers of $(w-a_h)$, while the constants $F^{n,g}_{k,m}[\mathcal{A}]$ are expressed in terms of a multinomial sum.

The conclusion is thus that the expression presented here is novel in its form. The usefulness of either expression will depend on the specific purposes for which it is required. An example for where the polynomial (\ref{Polynomial}) is more applicable is in the context of interpolation problems where the numerical values of the function and its derivatives at specific points are known, but not its general form. Expanding the derivatives in the constant $A_{n,g}$ as given by Eq. (\ref{A_ng2}) would ultimately result in expressions similar to the factors $F^{n,g}_{k,m}[\mathcal{A}]$, but this can be computationally expensive to produce and the resulting terms would not be grouped so naturally (for example, additional algebraic manipulations would be needed to obtain the factors (\ref{comp2}) which are useful for optimizing the computation of the polynomial). Thus, the expression (\ref{Polynomial}) can be more appropriate for explicit computations.

\section*{Acknowledgements}
Thanks to Stefan Nellen for proofreading the manuscript.

\section*{Financial disclosure}
Acknowledgements to the project DGAPA-UNAM IN103319 for financial support.

\section*{Conflict of interest}
The author declares no potential conflict of interests.

%\bibliography{references}

\begin{thebibliography}{10}

\bibitem{Richard}
R.~L. Burden and J.~D. Faires, {\em Numerical Analysis}.
\newblock Brooks/Cole, Cengage Learning, 9th~ed., 2010.

\bibitem{Philip}
P.~J. Davis, {\em Interpolation \& Approximation}.
\newblock Dover Publications, Inc., {D}over~ed., 1975.

\bibitem{Bulirsch}
J.~Stoer and R.~Bulirsch, {\em Introduction to Numerical Analysis}.
\newblock Springer-Verlag, 2nd~ed., 1993.

\bibitem{Hildebrand}
F.~B. Hildebrand, {\em Introduction to Numerical Analysis}.
\newblock Dover Publications, Inc., {D}over~ed., 1987.

\bibitem{Note_2Point}
K.~Kitahara, T.~Chiyonobu, and H.~Tsukamoto, ``A note on two point {T}aylor
  expansion,'' {\em International Journal of Pure and Applied Mathematics},
  vol.~75, no.~3, pp.~327--338, 2012.

\bibitem{NASA}
R.~H. Estes and E.~R. Lancaster, ``Two-point {T}aylor series expansions,''
  tech. rep., NASA Technical Reports Server, Dec 1966.

\bibitem{Lopez}
J.~L. López, E.~Pérez, and N.~M. Temme, ``Solving one-dimensional linear
  boundary value problems by multi-point {T}aylor polynomials. {A}pplications
  to special functions,'' {\em Monografías Matemáticas García de Galdeano},
  pp.~181--188, Sep 2010.

\bibitem{Diff_eq}
D.~S. Zézé {\em et~al.}, ``Multi-point {T}aylor series to solve differential
  equations,'' {\em Discrete and continuous dynamical systems series S},
  vol.~12, pp.~1791--1806, Oct 2019.

\bibitem{Image_comp}
G.~Franssens, M.~de~Maziere, D.~Fonteyn, and D.~Fussen, ``Image compression
  based on a multipoint {T}aylor series representation,'' in {\em Proceedings
  SIBGRAPI'98. International Symposium on Computer Graphics, Image Processing,
  and Vision (Cat. No.98EX237)}, pp.~174--184, 1998.

\bibitem{Gen_Leibniz}
R.~Majumdar, ``Generalization of {P}ascal's {R}ule and {L}eibniz's {R}ule for
  {D}ifferentiation,'' {\em Rose-Hulman Undergraduate Mathematics Journal},
  vol.~18, p.~9, Jul 2017.

\bibitem{franssens1999}
G.~Franssens, ``A new non-polynomial univariate interpolation formula of
  hermite type,'' {\em Advances in Computational Mathematics}, vol.~10,
  pp.~367--388, May 1999.

\bibitem{fine1962}
M.~Fine and J.~W. Head, ``A note on the convergence of multi-point {T}aylor's
  series,'' {\em Mathematical Proceedings of the Cambridge Philosophical
  Society}, vol.~58, pp.~548--550, Jul 1962.

\bibitem{2Point_N_Temme}
J.~L. López and N.~M. Temme, ``Two-point {T}aylor expansions of analytic
  functions,'' {\em Studies in Applied Mathematics}, vol.~109, no.~4,
  pp.~297--311, 2002.

\bibitem{MPoint_N_Temme}
J.~L. López and N.~M. Temme, ``Multi-point {T}aylor expansions of analytic
  functions,'' {\em Transactions of the American Mathematical Society},
  vol.~356, 11 2004.

\end{thebibliography}
%\bibliographystyle{ieeetr}

\section*{Appendix A: Proof that Eq. (\ref{Non_trivial}) is null.}\label{AppendixA}
Here it will be shown that the sum
\begin{equation}
\label{Sum}
\begin{split}
    \sum_{q=0}^{n+1}(-1)^{q}\binom{k+n+1-q}{n+1-q}\binom{k+1}{q}, \\
\end{split}
\end{equation}
is null. A renaming of the variables given by $n_j = n$ and $k_1 = q$ with respect to Eq. (\ref{Non_trivial}) was used for simplicity. For the calculation it is convenient to take the partial sums
\begin{equation}
\begin{split}
    S_p & \equiv \sum_{q=0}^{p+1}T_{n}(q), \quad \quad T_{n}(q) \equiv (-1)^{q}\binom{k+n+1-q}{n+1-q}\binom{k+1}{q}, \qquad p=0,...,n, \\
\end{split}
\end{equation}
where Eq. (\ref{Sum}) corresponds to the case $p=n$. Additional rearrangements of the factorials in the binomial coefficients allows one to express $T_{n}(q)$ as
\begin{equation}
\label{T_nq}
    T_n(q) = \frac{k+1}{n+1}(-1)^{q}\binom{k+n+1-q}{n}\binom{n+1}{q}.
\end{equation}
With this in mind, it will now be shown that
\begin{equation}
\label{S_T}
    S_p = -\frac{(p+2)(k+n-p)}{(n+1)(k+1)}T_n(p+2)
\end{equation}
holds for the partial sum $p=0,...,n-1$. Note that the second binomial coefficient of $T_n(q)$ in Eq. (\ref{T_nq}) is well defined only for $q=0,...,n+1$ and so in relation (\ref{S_T}) one cannot take $p=n$. This relation can be proven by induction. The base case is
\begin{equation}
\begin{split}
    S_0 & = \frac{k+1}{n+1}\bigg[\binom{k+n+1}{n}-\binom{k+n}{n}\binom{n+1}{1}\bigg] \\
    & = -\frac{1}{(n+1)}\frac{(k+n)!}{(k-1)!(n-1)!}
\end{split}
\end{equation}
which agrees with Eq. (\ref{S_T}) for $p=0$. For the inductive step it is now assumed that Eq. (\ref{S_T}) holds for a given $p = 0,...,n-2$. This implies that
\begin{equation}
\begin{split}
    S_{p+1} & = S_p+T_n(p+2) \\
    & = \bigg[-\frac{(p+2)(k+n-p)}{(n+1)(k+1)}+1\bigg]T_n(p+2) \\
    & = \frac{(k-1-p)(n-1-p)}{(n+1)^2}(-1)^{p+2}\binom{k+n-1-p}{n}\binom{n+1}{p+2} \\
    & = \frac{(p+3)(k+n-1-p)}{(n+1)^2}(-1)^{p+2}\binom{k+n-2-p}{n}\binom{n+1}{p+3} \\
    & = -\frac{(p+3)(k+n-1-p)}{(n+1)(k+1)}T_n(p+3)
\end{split}
\end{equation}
and thus induction is fulfilled: Eq. (\ref{S_T}) holds for $p=0,...,n-1$. With this relation the complete sum can be now directly computed:
\begin{equation}
    S_n = S_{n-1}+T_n(n+1) = \bigg[-\frac{(n+1)(k+1)}{(n+1)(k+1)}+1\bigg]T_n(n+1) = 0.
\end{equation}
The sum vanishes, as required.
$\blacksquare$

\end{document}